%
%
%
%


\magnification=1200
\pretolerance=500 \tolerance=1000 \brokenpenalty=5000
\hsize=12.5cm   
\vsize=19cm
\hoffset=0.4cm
\voffset=1cm
\parskip3pt plus 1pt
\parindent=0.4cm
\let\sl=\it
\def\\{\hfil\break}


\font\seventeenbf=cmbx10 at 17.28pt

\font\twelvebf=cmbx10 at 12pt
\font\eightbf=cmbx8
\font\sixbf=cmbx6

\font\eighti=cmmi8
\font\sixi=cmmi6

\font\eightrm=cmr8
\font\sixrm=cmr6

\font\eightsy=cmsy8
\font\sixsy=cmsy6

\font\eightit=cmti8
\font\eighttt=cmtt8
\font\eightsl=cmsl8

\font\seventeenbsy=cmbsy10 at 17.28pt

\font\twelvebsy=cmbsy10 at 12pt
\font\tenbsy=cmbsy10
\font\eightbsy=cmbsy8
\font\sevenbsy=cmbsy7
\font\sixbsy=cmbsy6
\font\fivebsy=cmbsy5

\font\tenmsa=msam10

\font\sevenmsa=msam7
\font\fivemsa=msam5
\newfam\msafam
  \textfont\msafam=\tenmsa
  \scriptfont\msafam=\sevenmsa
  \scriptscriptfont\msafam=\fivemsa

\font\tenmsb=msbm10
\font\eightmsb=msbm8
\font\sevenmsb=msbm7
\font\fivemsb=msbm5
\newfam\msbfam
  \textfont\msbfam=\tenmsb
  \scriptfont\msbfam=\sevenmsb
  \scriptscriptfont\msbfam=\fivemsb
\def\Bbb{\fam\msbfam\tenmsb}

\font\tenCal=eusm10
\font\sevenCal=eusm7
\font\fiveCal=eusm5
\newfam\Calfam
  \textfont\Calfam=\tenCal
  \scriptfont\Calfam=\sevenCal
  \scriptscriptfont\Calfam=\fiveCal
\def\Cal{\fam\Calfam\tenCal}

\font\teneuf=eusm10
\font\teneuf=eufm10
\font\seveneuf=eufm7
\font\fiveeuf=eufm5
\newfam\euffam
  \textfont\euffam=\teneuf
  \scriptfont\euffam=\seveneuf
  \scriptscriptfont\euffam=\fiveeuf

\font\seventeenbfit=cmmib10 at 17.28pt

\font\twelvebfit=cmmib10 at 12pt
\font\tenbfit=cmmib10
\font\eightbfit=cmmib8
\font\sevenbfit=cmmib7
\font\sixbfit=cmmib6
\font\fivebfit=cmmib5
\newfam\bfitfam
  \textfont\bfitfam=\tenbfit
  \scriptfont\bfitfam=\sevenbfit
  \scriptscriptfont\bfitfam=\fivebfit


\catcode`\@=11
\def\eightpoint{%
  \textfont0=\eightrm \scriptfont0=\sixrm \scriptscriptfont0=\fiverm
  \def\rm{\fam\z@\eightrm}%
  \textfont1=\eighti \scriptfont1=\sixi \scriptscriptfont1=\fivei
  \def\oldstyle{\fam\@ne\eighti}%
  \textfont2=\eightsy \scriptfont2=\sixsy \scriptscriptfont2=\fivesy
  \textfont\itfam=\eightit
  \def\it{\fam\itfam\eightit}%
  \textfont\slfam=\eightsl
  \def\sl{\fam\slfam\eightsl}%
  \textfont\bffam=\eightbf \scriptfont\bffam=\sixbf
  \scriptscriptfont\bffam=\fivebf
  \def\bf{\fam\bffam\eightbf}%
  \textfont\ttfam=\eighttt
  \def\tt{\fam\ttfam\eighttt}%
  \textfont\msbfam=\eightmsb
  \def\Bbb{\fam\msbfam\eightmsb}%
  \abovedisplayskip=9pt plus 2pt minus 6pt
  \abovedisplayshortskip=0pt plus 2pt
  \belowdisplayskip=9pt plus 2pt minus 6pt
  \belowdisplayshortskip=5pt plus 2pt minus 3pt
  \smallskipamount=2pt plus 1pt minus 1pt
  \medskipamount=4pt plus 2pt minus 1pt
  \bigskipamount=9pt plus 3pt minus 3pt
  \normalbaselineskip=9pt
  \setbox\strutbox=\hbox{\vrule height7pt depth2pt width0pt}%
  \let\bigf@ntpc=\eightrm \let\smallf@ntpc=\sixrm
  \normalbaselines\rm}
\catcode`\@=12

\def\eightpointbf{%
 \textfont0=\eightbf   \scriptfont0=\sixbf   \scriptscriptfont0=\fivebf
 \textfont1=\eightbfit \scriptfont1=\sixbfit \scriptscriptfont1=\fivebfit
 \textfont2=\eightbsy  \scriptfont2=\sixbsy  \scriptscriptfont2=\fivebsy
 \eightbf
 \baselineskip=10pt}

\def\tenpointbf{%
 \textfont0=\tenbf   \scriptfont0=\sevenbf   \scriptscriptfont0=\fivebf
 \textfont1=\tenbfit \scriptfont1=\sevenbfit \scriptscriptfont1=\fivebfit
 \textfont2=\tenbsy  \scriptfont2=\sevenbsy  \scriptscriptfont2=\fivebsy
 \tenbf}
        
\def\twelvepointbf{%
 \textfont0=\twelvebf   \scriptfont0=\eightbf   \scriptscriptfont0=\sixbf
 \textfont1=\twelvebfit \scriptfont1=\eightbfit \scriptscriptfont1=\sixbfit
 \textfont2=\twelvebsy  \scriptfont2=\eightbsy  \scriptscriptfont2=\sixbsy
 \twelvebf
 \baselineskip=14.4pt}

\def\seventeenpointbf{%
 \textfont0=\seventeenbf  \scriptfont0=\twelvebf  \scriptscriptfont0=\eightbf
 \textfont1=\seventeenbfit\scriptfont1=\twelvebfit\scriptscriptfont1=\eightbfit
 \textfont2=\seventeenbsy \scriptfont2=\twelvebsy \scriptscriptfont2=\eightbsy
 \seventeenbf
 \baselineskip=20.736pt}
 
       
\newbox\titlebox   \setbox\titlebox\hbox{\hfil}
\newbox\sectionbox \setbox\sectionbox\hbox{\hfil}
\def\folio{\ifnum\pageno=1 \hfil \else \ifodd\pageno
           \hfil {\eightpoint\copy\sectionbox\kern8mm\number\pageno}\else
           {\eightpoint\number\pageno\kern8mm\copy\titlebox}\hfil \fi\fi}
\footline={\hfil}
\headline={\folio}

\def\titlerunning#1{\setbox\titlebox\hbox{\eightpoint #1}}
\def\title#1{\noindent\hfil$\smash{\hbox{\seventeenpointbf #1}}$\hfil
             \titlerunning{#1}\medskip}

\newcount\numbersection \numbersection=-1
\def\sectionrunning#1{\setbox\sectionbox\hbox{\eightpoint #1}}
\def\section#1{%
  \par\vskip0.666cm\penalty -100
  \vbox{\baselineskip=14.4pt\noindent{{\twelvepointbf #1}}}
  \vskip2pt
  \penalty 500
  \advance\numbersection by 1
  \sectionrunning{#1}}

\def\subsection#1|{%
  \par\vskip0.5cm\penalty -100
  \vbox{\noindent{{\tenpointbf #1}}}
  \vskip1pt
  \penalty 500}

\newcount\numberindex \numberindex=0  
\def\index#1#2{%
  \advance\numberindex by 1
  \immediate\write1{\string\def \string\IND #1%
     \romannumeral\numberindex \string{%
     \noexpand#2 \string\dotfill \space \string\S \number\numbersection, 
     p.\string\ \space\number\pageno \string}}}

\newdimen\itemindent \itemindent=\parindent

\def\item#1{\par\noindent\hangindent\itemindent%
            \rlap{#1}\kern\itemindent\ignorespaces}
\def\itemitem#1{\par\noindent\hangindent2\itemindent%
            \kern\itemindent\rlap{#1}\kern\itemindent\ignorespaces}
\def\itemitemitem#1{\par\noindent\hangindent3\itemindent%
            \kern2\itemindent\rlap{#1}\kern\itemindent\ignorespaces}

\long\def\claim#1|#2\endclaim{\par\vskip 5pt\noindent 
{\tenpointbf #1.}\ {\sl #2}\par\vskip 5pt}

\def\proof{\noindent{\sl Proof}}

\def\today{\ifcase\month\or
January\or February\or March\or April\or May\or June\or July\or August\or
September\or October\or November\or December\fi \space\number\day,
\number\year}

\catcode`\@=11
\newcount\@tempcnta \newcount\@tempcntb 
\def\timeofday{{%
\@tempcnta=\time \divide\@tempcnta by 60 \@tempcntb=\@tempcnta
\multiply\@tempcntb by -60 \advance\@tempcntb by \time
\ifnum\@tempcntb > 9 \number\@tempcnta:\number\@tempcntb
  \else\number\@tempcnta:0\number\@tempcntb\fi}}
\catcode`\@=12

\def\bibitem#1&#2&#3&#4&%
{\hangindent=1.66cm\hangafter=1
\noindent\rlap{\hbox{\eightpointbf #1}}\kern1.66cm{\rm #2}{\sl #3}{\rm #4.}} 


\def\bC{{\Bbb C}}

\def\bN{{\Bbb N}}

\def\bQ{{\Bbb Q}}
\def\bR{{\Bbb R}}

\def\bZ{{\Bbb Z}}


\def\cC{{\Cal C}}
\def\cE{{\Cal E}}
\def\cF{{\Cal F}}
\def\cG{{\Cal G}}

\def\cO{{\Cal O}}
\def\cI{{\Cal I}}
\def\cJ{{\Cal J}}
\def\cU{{\Cal U}}


\def\square{{\hfill \hbox{
\vrule height 1.453ex  width 0.093ex  depth 0ex
\vrule height 1.5ex  width 1.3ex  depth -1.407ex\kern-0.1ex
\vrule height 1.453ex  width 0.093ex  depth 0ex\kern-1.35ex
\vrule height 0.093ex  width 1.3ex  depth 0ex}}}
\def\bigsquare{{\kern-0.3ex\hbox{
\vrule height 1.7ex  width 0.093ex  depth 0ex\kern-0.093ex
\vrule height 1.8ex  width 1.7ex  depth -1.707ex\kern-0.093ex
\vrule height 1.7ex  width 0.093ex  depth 0ex\kern-1.65ex
\vrule height 0.093ex  width 1.6ex  depth 0ex}\kern0.3ex}}
\def\qed{\phantom{$\quad$}\hfill$\square$\medskip}
\def\hexnbr#1{\ifnum#1<10 \number#1\else
 \ifnum#1=10 A\else\ifnum#1=11 B\else\ifnum#1=12 C\else
 \ifnum#1=13 D\else\ifnum#1=14 E\else\ifnum#1=15 F\fi\fi\fi\fi\fi\fi\fi}
\def\msatype{\hexnbr\msafam}
\def\msbtype{\hexnbr\msbfam}
\mathchardef\restriction="3\msatype16   
\mathchardef\compact="3\msatype62
\mathchardef\smallsetminus="2\msbtype72   \let\ssm\smallsetminus
\mathchardef\subsetneq="3\msbtype28
\mathchardef\supsetneq="3\msbtype29
\mathchardef\leqslant="3\msatype36   \let\le\leqslant
\mathchardef\geqslant="3\msatype3E   \let\ge\geqslant
\mathchardef\ltimes="2\msbtype6E
\mathchardef\rtimes="2\msbtype6F


\let\text=\hbox
\def\build#1|#2|#3|{\mathrel{\mathop{\null#1}\limits^{#2}_{#3}}}

\def\Ch{\mathop{\rm Ch}\nolimits}
\def\Todd{\mathop{\rm Todd}\nolimits}

\def\Pic{\mathop{\rm Pic}\nolimits}

\def\Supp{\mathop{\rm Supp}\nolimits}

\def\Vol{\mathop{\rm Vol}\nolimits}

\def\Sing{\mathop{\rm sing}\nolimits}

\def\ddbar{{\partial\overline\partial}}

\def\red{{\rm red}}

\def\BC{{\rm BC}}
\def\DR{{\rm DR}}
\def\NS{{\rm NS}}
\def\DNS{{\rm DNS}}
\def\tr{{\rm tr}}

\long\def\InsertFig#1 #2 #3 #4\EndFig{\par
\hbox{\hskip #1mm$\vbox to#2mm{\vfil\special{" 
(/home/demailly/psinputs/grlib.ps) run
#3}}#4$}}
\long\def\LabelTeX#1 #2 #3\ELTX{\rlap{\kern#1mm\raise#2mm\hbox{#3}}}


\title{Holomorphic Morse inequalities}
\smallskip
\title{ and asymptotic cohomology groups:}
\smallskip
\title{a tribute to Bernhard Riemann}
\titlerunning{Asymptotic cohomological inequalities: a tribute to Riemann}
\medskip
\centerline{\twelvebf Jean-Pierre Demailly}
\medskip
\centerline{Universit\'e de Grenoble I, D\'epartement de Math\'ematiques}
\centerline{Institut Fourier, 38402 Saint-Martin d'H\`eres, France}
\centerline{{\it e-mail\/}: {\tt demailly@fourier.ujf-grenoble.fr}}

\smallskip
\vskip35pt

\noindent
{\bf Abstract.}
The goal of this note is to present the potential relationships 
between  certain Monge-Amp\`ere integrals appearing in holomorphic Morse
inequa\-lities, and asymptotic cohomology estimates for tensor
powers of line bundles, as recently introduced by algebraic geometers.
The expected most general statements are still conjectural and
owe a debt to Riemann's pioneering work, which led 
to the concept of Hilbert polynomials and to the Hirzebruch-Riemann-Roch
formula during the XX-th century.
\medskip

\noindent
{\bf R\'esum\'e.}
Le but de cette note est de pr\'esenter les relations potentielles qui
doivent exister entre certaines int\'egrales de Monge-Amp\`ere et les
estimations asymptotiques de cohomologie introduites r\'ecemment par
les g\'eom\`etres alg\'ebristes. Les \'enonc\'es les plus g\'en\'eraux
esp\'er\'es sont encore conjecturaux, et ne peuvent \^etre formul\'es
sans faire r\'ef\'erence aux travaux pionniers de Riemann, qui ont
conduit au concept de polyn\^ome de Hilbert et \`a la formule de
Hirzebruch-Riemann-Roch au cours du XXe si\`ecle.
\medskip\noindent
{\bf Mathematics Subject Classification 2010.} 32F07, 14B05, 14C17
\medskip\noindent
{\bf Key words.}
Holomorphic Morse inequalities, Monge-Amp\`ere integrals, Dolbeault
cohomology, asymptotic cohomology groups, Riemann-Roch formula,
hermitian metrics, Chern curvature tensor, plurisuharmonic
approximation.
\medskip\noindent
{\bf Mots-cl\'es.}
In\'egalit\'es de Morse holomorphes, int\'egrales de Monge-Amp\`ere, 
cohomologie de Dolbeault, formule de Riemann-Roch, m\'etriques
hermitiennes, courbure de Chern, approximation plurisousharmonique.
\vfill\eject

\section{1. Main results}

Throughout the paper, $X$ will denote a compact complex manifold and
\hbox{$n=\dim_\bC X$} its complex dimension. Hirzebruch's Riemann-Roch 
formula [Hir54, Hir56] expresses the Euler characteristic 
$$
\chi(X,\cF)=\sum_{j=0}^n(-1)^j h^j(X,\cF)\leqno(1.1)
$$
of any coherent analytic sheaf $\cF$ over $X$ as an explicit topological
invariant computed by the integral
$$
\int_X\Ch(\cF)\Todd(T_X)\leqno(1.2)
$$
in terms of the Chern character of $\cF$ and the Todd class of $T_X$. 
In the special case where $\cF=\cO(L^{\otimes k})$ is the $k$-th tensor power
of a holomorphic line bundle, the formula produces the Hilbert polynomial
$$
\chi(X,L^{\otimes k})=P(k)\leqno(1.3)
$$
which is a polynomial of degree $n$ which leading term ${k^n\over n!}c_1(L)^n$.
Moreover, if $h$ is a hermitian metric on $L$ and $\Theta_{L,h}={i\over 2\pi}
D_h^2$ is the $(1,1)$ Chern curvature tensor of~$(L,h)$, the top Chern class
intersection number is given by
$$
c_1(L)^n={k^n\over n!}\int_X\Theta_{L,h}^n.\leqno(1.4)
$$
However, the Hilbert polynomial just gives access to the alternate sum
of dimensions; in order to estimate the growth of the
individual cohomology groups, it is interesting to consider
appropriate ``asymptotic cohomology functions''. We mostly follow 
here notation and concepts 
introduced by A.~K\"uronya [Kur06, FKL07], which extend in 
a natural way the concept of volume of a line bundle (cf.\ [DEL00], 
[Bou02], [Laz04]).

However, as we insist here on being able to deal with arbitrary
compact complex manifolds, we are led to introduce appropriate 
variants of the original definitions.  Recall that the Bott-Chern cohomology
group $H^{p,q}_\BC(X,\bC)$ is the quotient of $d$-closed $(p,q)$-forms by
$\smash{\ddbar}$-exact $(p,q)$-forms. Then
$\bigoplus_{p,q}H^{p,q}_\BC(X,\bC)$ is a bigraded finite dimensional
algebra. When $X$ is K\"ahler, this algebra is isomorphic to the usual 
Hodge-De Rham cohomology algebra by the well-known $\ddbar$-lemma, 
but in general we only have canonical morphisms
$$
H^{p,q}_\BC(X,\bC)\to H^{p,q}(X,\bC),\qquad 
\bigoplus_{p+q=k} H^{p,q}_\BC(X,\bC)\to H^k_\DR(X,\bC)
$$
to Dolbeault and de Rham cohomology groups, which need not be isomorphisms.
The Bott-Chern cohomology algebra carries a natural
conjugation, and we can thus look at real elements
$H^{p,p}_\BC(X,\bR)$ of type $(p,p)$. The first Chern class of a
holomorphic line bundle $L\to X$ yields a well defined 
Bott-Chern class \hbox{$c_1(L)\in H^{1,1}_\BC(X,\bR)$} and conversely,
by a well known lemma due to A.~Weil, such classes correspond 
to a holomorphic line bundle if and only if they are integral,
i.e.\ in the image of $H^2(X,\bZ)\to H^2(X,\bR)$ under
the canonical morphism $H^{1,1}_\BC(X,\bR)\to H^2(X,\bR)$.

We consider the real Neron-Severi subspace $\NS_\bR(X)\subset
H^{1,1}_\BC(X,\bR)$ generated by real combinations of all
Chern classes $c_1(L)$ (i.e., from what we said, by integral $(1,1)$
classes). Given a cohomology
class $\alpha\in \NS_\bR(X)$ there is always a sequence of $\bQ$-line
bundles $\smash{1\over
  k_\nu}L_{k_\nu}\in\Pic_\bQ(X)=\Pic(X)\otimes_\bZ\bQ$ such that
$\smash{1\over k_\nu}c_1(L_{k_\nu})$ converges to $\alpha$ in
$\NS_\bR(X)$.  We will simply write ${1\over k}c_1(L)\to \alpha$ to
express the fact that ${1\over k}c_1(L)$ is close to $\alpha$ in the
finite dimensional vector space $\NS_\bR(X)\subset
H^{1,1}_\BC(X,\bR)$, with its natural Hausdorff topology.

\claim (1.5) Definition|Let $X$ be a compact complex manifold. One
defines the asymptotic $($analytic$)$ $q$-cohomology function on 
$\NS_\bR(X)$ to be
$$
\eqalign{
\widehat h^q(X,\alpha)&=
\limsup_{k\to+\infty, {1\over k}c_1(L)\to \alpha}~{n!\over k^n}h^q(X,L)\cr
&=\inf_{\varepsilon>0,\,k_0>0}~~
\sup_{k\ge k_0,\Vert{1\over k}c_1(L)-\alpha\Vert\le\varepsilon}
~{n!\over k^n}h^q(X,L).\cr}
$$
where the pair $(k,L)$ runs over $\bN^*\times\Pic(X)$.
\endclaim

From the very definition, $\widehat h^q$ is an upper semi-continuous 
function on $NS_\bR(X)$ and it is positively homogeneous of degree $n$, namely 
$$
\widehat h^q(X,\lambda\alpha)=\lambda^n\widehat h^q(X,\alpha)\leqno(1.6)
$$
for all $\alpha\in\NS_\bR(X)$ and all $\lambda\ge 0$. In the general
case of compact complex manifolds, the fact that $\widehat
h^q(X,\alpha)$ is finite follows from spectral theory estimates for
the complex Laplace-Beltrami operators (this will become quite
clear from the discussions below).

For a line bundle $L$, we simply denote $\widehat h^q(X,L)=
\widehat h^q(X,c_1(L))$. In this case, things become a little bit simpler,
and especially, for $q=0$, one recovers the usual concept of volume
of a line bundle.

\claim (1.7) Proposition|If $X$ is projective algebraic or $q=0$, then
$$
\widehat h^q(X,L)=
\limsup_{k\to+\infty}~~{n!\over k^n}h^q(X,L^{\otimes k})
= \lim_{k\to+\infty}~~{n!\over k^n}h^q(X,L^{\otimes k}).
$$
Moreover, in these cases, the map $\alpha\mapsto\widehat h^q(X,\alpha)$ is
$($locally$)$ Lipschitz continuous on $\NS_\bR(X)$.
\endclaim

The proof is derived from arguments quite similar to those already 
developed in [Kur05]. Actually,
let us introduce $\DNS_\bR(X)\subset \NS_\bR(X)$ to be the subspace generated
by classes of integral divisors $D$ on $X$ (``divisorial Neron-Severi group'').
If $X$ is projective algebraic
then $\DNS_\bR(X)=\NS_\bR(X)$, but the inclusion can be strict in general
(e.g.\ on complex 2-tori which only have indefinite 
integral $(1,1)$-classes, cf. [BL04]). If $D=\sum p_jD_j$ is
an integral divisor, we define its norm to be $\Vert D\Vert=\sum |p_j|
\Vol_\omega(D_j)$, where the volume of an irreducible divisor is computed
by means of a given hermitian metric $\omega$ on~$X$; in other words, this
is precisely the mass of the current of integration $[D]$ with 
respect to~$\omega$.
Clearly, since $X$ is compact, we get equivalent norms for all
choices of hermitian metrics $\omega$ on~$X$. We can also use $\omega$
to fix a normalized metric on $H^{1,1}_\BC(X,\bR)$. Elementary properties
of potential theory show that $\Vert c_1(\cO(D))\Vert\le C\Vert D\Vert$
for some constant $C>0$ (but the converse inequality is of course wrong
in most cases). Proposition 1.7 
is a simple consequence of the following more precise cohomology estimates 
which will be proved in section~2. 

\claim (1.8) Theorem|Let $X$ be a compact complex manifold. Fix a finitely 
generated subgroup $\Gamma$ of the group of $\bZ$-divisors on~$X$. Then
there are constants $C$, $C'$ depending only on~$X$, its hermitian metric 
$\omega$ and the subgroup $\Gamma$, satisfying the following properties.
\smallskip
{\itemindent=6.5mm
\item{\rm(a)} Let $L$ and $L'=L\otimes\cO(D)$ be holomorphic line bundles
on $X$, where $D\in\Gamma$ is an integral divisor. Then 
$$
\big| h^q(X,L')-h^q(X,L)\big|
\le C(\Vert c_1(L)\Vert+\Vert D\Vert)^{n-1}\Vert D\Vert.
$$
\item{\rm(b)} On the subspace $\DNS_\bR(X)$, the asymptotic $q$-cohomology
function $\widehat h^q$ satisfies a global estimate
$$
\big|\widehat h^q(X,\beta)-\widehat h^q(X,\alpha)\big|\le C'
(\Vert \alpha\Vert+\Vert\beta\Vert)^{n-1}\Vert \beta-\alpha\Vert.
$$}
\vskip-\baselineskip\noindent
In particular $($without any further assumption on $X)$, $\widehat h^q$
is locally Lipschitz continuous on $\DNS_\bR(X)$.
\endclaim

Our ambition is to extend the function $\widehat h^q$ in a
natural way to the full cohomology group $H^{1,1}_\BC(X,\bR)$. The main
trouble, already when $X$ is projective algebraic, is that the Picard
number $\rho(X)=\dim_\bR\NS_\bR(X)$ may be much smaller than
$\dim_\bR H^{1,1}_\BC(X,\bR)$, namely, there can be rather few integral classes
of type~$(1,1)$ on~$X$. It is well known for instance that $\rho(X)=0$ for
a generic complex torus a dimension $n\ge 2$, while
$\dim_\bR H^{1,1}_\BC(X,\bR)=n^2$. However, if we look at the natural 
morphism
$$
H^{1,1}_\BC(X,\bR)\to H^2_\DR(X,\bR)\simeq H^2(X,\bR)
$$
to de Rham cohomology, then $H^2(X,\bQ)$ is dense in $H^2(X,\bR)$. Therefore,
given a class $\alpha\in H^{1,1}_\BC(X,\bR)$ and a smooth $d$-closed
$(1,1)$-form $u$ in $\alpha$,
we can find an infinite sequence ${1\over k}L_k$ ($k\in S\subset\bN$) of 
topological $\bQ$-line bundles, equipped with hermitian metrics $h_k$ and
compatible connections $\nabla_k$ such that the curvature forms 
${1\over k}\Theta_{\nabla_k}$ converge to $u$. By using Kronecker's approximation
with respect to the integral lattice $H^2(X,\bZ)/{\rm torsion}
\subset H^2(X,\bR)$, we can even achieve a fast diophantine approximation 
$$\Vert  \Theta_{\nabla_k}-ku\Vert\le Ck^{-1/b_2}\leqno(1.9)$$
for a suitable infinite subset $k\in S\subset\bN$ of multipliers. Then
in particular
$$\Vert\Theta_{\nabla_k}^{0,2}\Vert=\Vert\Theta_{\nabla_k}^{0,2}-u^{0,2}\Vert\le
Ck^{-1/b_2},$$
and we see that $(L_k,h_k,\nabla_k)$ is a $C^\infty$ hermitian 
line bundle which is extremely close to being holomorphic, since
$(\nabla_k^{0,1})^2=\Theta_{\nabla_k}^{0,2}$ is very small. We introduce
the complex Laplace-Beltrami operator
$$
\overline{\bigsquare}_k=(\nabla_k^{0,1})(\nabla_k^{0,1})^*+(\nabla_k^{0,1})^*(\nabla_k^{0,1})
$$
and look at its eigenspaces in $L^2(X,\Lambda^{0,q}T^\star X\otimes L_k)$
with the metric induced by $\omega$ on $X$ and $h_k$ on $L_k$. In the 
holomorphic case, Hodge theory tells us that the $0$-eigenspace is
isomorphic to $H^q(X,\cO(L_k))$, but in the ``almost holomorphic case'' 
the $0$-eigenvalues deviate from $0$, essentially by a
shift of the order of magnitude of 
$\Vert\Theta_{\nabla_k}^{0,2}\Vert\sim k^{-1/b_2}$
(see [Lae02], chapter~4). It is thus natural to introduce
in this case

\claim (1.10) Definition|Let $X$ be a compact complex manifold and
$\alpha\in H^{1,1}_\BC(X,\bR)$ an arbitrary Bott-Chern $(1,1)$-class.
We define the ``transcendental'' asymptotic $q$-cohomology function
to be
$$
\widehat h^q_\tr(X,\alpha)=\inf_{u\in\alpha}~~
\limsup_{\varepsilon\to 0,\,k\to+\infty,\,L_k,\,h_k,\,\nabla_k,
{1\over k}\Theta_{\nabla_k}\to u}~~{n!\over k^n}
N(\overline{\bigsquare}_k,k\varepsilon)
$$
where the $\limsup$ runs over all $5$-tuples $(\varepsilon,k,L_k,h_k,\nabla_k)$,
and where $N(\overline{\bigsquare}_k,k\varepsilon)$ denotes the sum of
dimensions of all eigenspaces of eigenvalues at most equal to $k\varepsilon$ 
for the Laplace-Beltrami operator
$\overline{\bigsquare}_k$ associated with $(L_k,h_k,\nabla_k)$ and the
base hermitian metric~$\omega$.
\endclaim

The word ``transcendental'' refers here to the fact that we deal with classes
$\alpha$ of type $(1,1)$ which are not algebraic or even analytic. Of course,
in the definition, we could have restricted the limsup to families satisfying 
a better approximation property $\Vert{1\over k}\Theta_{\nabla_k}-u\Vert\le C
k^{-1-1/b_2}$ for some large constant $C$ (this would lead a priori to
a smaller limsup, but there is enough stability in the parameter
dependence of the spectrum for making such a change irrelevant). 
The minimax principle easily shows that definition 1.10 does not 
depend on $\omega$, as the eigenvalues
are at most multiplied or divided by constants under a change of base metric.
When $\alpha\in\NS_\bR(X)$, by restricting our families 
$\{(\varepsilon,k,L_k,h_k,\nabla_k)\}$ to the case
of holomorphic line bundles only, we get the obvious inequality
$$
\widehat h^q(X,\alpha)\le \widehat h^q_\tr(X,\alpha),\qquad
\forall\alpha\in\NS_\bR(X).\leqno(1.11)
$$
It is natural to raise the question whether this is always an equality. 
Hopefully, the calculation of the quantities
 $\limsup~{n!\over k^n}
N(\overline{\bigsquare}_k,k\varepsilon)$ is a problem of spectral theory 
which is completely understood since a long time. In fact, as
a consequence of the techniques of [Dem85, Dem91, Lae02], one gets

\claim (1.12) Theorem|With the above notations and assumptions, one has
$$
\limsup_{\varepsilon\to 0,\,k\to+\infty,\,L_k,\,h_k,\,\nabla_k,
{1\over k}\Theta_{\nabla_k}\to u}~~{n!\over k^n}
N(\overline{\bigsquare}_k,k\varepsilon)=
\int_{X(u,q)}(-1)^qu^n,
$$
where $X(u,q)$ is the open set of points $x\in X$ where $u(x)$ has signature
$(n-q,q)$. Therefore
$$
\widehat h^q_\tr(X,\alpha)=\inf_{u\in\alpha}\int_{X(u,q)}(-1)^qu^n\qquad
(\hbox{$u$ smooth}).
$$
\endclaim

The first equality follows mainly from Theorems~2.16 and 3.14 of [Dem85], 
which even yield explicitly the limit for any given $\varepsilon$ outside 
a countable set (the limit as $\varepsilon\to 0$ is then obtained from the 
calculations of page 224 after Cor.~4.3). One has to observe, in the case 
of sequences of ``almost holomorphic line bundles'' considered here, that 
the perturbation indeed goes to $0$, and also that 
all constants involved in the calculations of [Dem85] are uniformly 
bounded; see [Dem91] and [Lae02] for more details on this. Therefore,
we can reformulate more explicitly our previous question in the following 
terms.

\claim (1.13) Question|For every $\alpha\in\NS_\bR(X)$, is it true that
$$
\widehat h^q(X,\alpha)=
\inf_{u\in\alpha}\int_{X(u,q)}(-1)^qu^n\qquad
(\hbox{$u$ smooth})~?
$$
{\rm (Note: it is known, from the holomorphic Morse inequalities proved 
in [Dem85], that the inequality $\le$ always holds true).}
\endclaim

In general, equality (1.13) seems rather hard to prove. In some sense,
this would be an asymptotic converse of the Andreotti-Grauert theorem
[AG62]~: under a suitable $q$-convexity assumption, the latter asserts
the vanishing of related cohomology groups in degree~$q$; here,
conversely, assuming a known growth of these groups in degree~$q$, we
expect to be able to say something about the $q$-index sets of
suitable hermitian metrics on the line bundles under consideration.

For degree $q=0$, however, we deal with sections rather than with cohomology
classes, and complex pluripotential theory makes things much easier. In the
case $q=0$, there are for instance some well known methods to compute
the volume $\Vol(\alpha)$ of a transcendental class $\alpha\in
H^{1,1}_\BC(X,\bR)$. 

\claim (1.14) Definition|Let $X$ be a compact complex $n$-fold.
We denote by 
\smallskip
{\itemindent=6.5mm
\item{\rm(a)} $\cE_X\subset H^{1,1}_\BC(X,\bR)$ the 
{\rm pseudoeffective cone} of~$X$, namely the cone of 
classes of closed positive $(1,1)$-currents; this is a closed
convex cone$\,;$
\smallskip
\item{\rm(b)} $\cE_X^+\subset\cE_X$ the cone
consisting of classes of {\rm K\"ahler currents}, i.e.\ positive
currents which admit a positive lower bound $T\ge\varepsilon\omega$ where
$\omega$ is a smooth positive $(1,1)$-form on $X$ and $\varepsilon>0\,;$
this is an open convex cone.\vskip0pt}
\endclaim

Given a class $\alpha\in H^{1,1}_\BC(X,\bR)$, we set $\Vol(\alpha)=0$ if 
$\alpha\notin\cE^+_X$. Otherwise, if $\alpha\in\cE^+_X$, the main 
approximation theorem 
of [Dem92] shows that the class
$\alpha$ contains K\"ahler currents $T$ with analytic singularities, i.e.\
such that their local potentials $\varphi$
of $T$ have singularities of the form $\varphi={1\over k}\log|\sum_j|g_{j,k}|^2$
mod $C^\infty$, for suitable local holomorphic functions $(g_{j,k})$.
Then there exists a blow-up $\mu:\smash{\widetilde X}\to X$ of $X$ such that
$\mu^*T=[E]+\beta$, where $E$ is a divisor supported on 
$\mu^{-1}(\bigcap g_{j,k}^{-1}(0))$ and $\beta$ a smooth closed
positive $(1,1)$-form on~$\smash{\widetilde X}$ (cf.\ [BDPP04]). 
One can define
$$
\leqalignno{
\Vol(T)&=\int_{X\ssm\Sing(T)}T^n=\int_{\widetilde X}\beta^n,&(1.15)\cr
\Vol(\alpha)&=\sup_{\varepsilon\omega\le T\in \alpha}\Vol(T),&(1.16)\cr}
$$
where the supremum is taken over all K\"ahler currents with analytic 
singularities in the class~$\alpha$. By definition, the volume
function is identically zero unless $X$ carries K\"ahler currents, and by
[DP04] the latter property is equivalent to $X$ being in the {\it Fujiki 
class}~$\cC$ of manifolds bimeromorphic to K\"ahler. The results of 
S.~Boucksom [Bou02]  yield:

\claim (1.17) Theorem {\rm([Bou02])}|If $X$ is compact complex manifold and 
$\alpha\in\NS_\bR(X)$, then
$$\Vol(\alpha)=\widehat h^0(X,\alpha).$$
In other words, the growth of
sections of multiples of a line bundle $L$ can be calculated as the sup of 
volumes of  K\"ahler currents $T\in c_1(L)$ as defined above.
\endclaim

In section 3, we use the results of [BDPP04] and [BD09] to derive a proof of
the following theorem, which, in combination with Boucksom's theorem, 
yields a positive answer to question~(1.13) when $q=0$ and $X$ is a 
projective surface.

\claim (1.18) Theorem|Let $(X,\omega)$ be a compact complex $n$-fold.
Then for every class $\alpha\in H^{1,1}(X,\bR)$ we have
{\itemindent=6.5mm
\item{\rm(a)}
$\displaystyle
\Vol(\alpha)\le \widehat h^0_\tr(X,\alpha)=\inf_{u\in\alpha}\int_{X(u,0)}u^n$,
\smallskip
\item{} where the infimum runs over all smooth closed $(1,1)$-forms $u$
contained in the class~$\alpha$.
\smallskip
\item{\rm (b)} Equality holds if $X$ is a projective surface and 
$\alpha\in\NS_\bR(X)$.
}
\endclaim

It would be interesting to knwow whether equality always holds without
restrictions on $X$ or on~$\alpha$. In~the general setting of compact
complex manifolds, we also hope for the following ``transcendental''
case of holomorphic Morse inequalities.

\claim (1.19) Conjecture|Let $X$ be a compact complex $n$-fold and
$\alpha$ an arbitrary cohomology class in~$H^{1,1}_\BC(X,\bR)$. Then
$$
\Vol(\alpha)\ge\sup_{u\in\alpha}\int_{X(u,0)\cup X(u,1)}u^n.
$$
In particular, if the right hand side is positive, then $\alpha$
contains a K\"ahler current and $X$ must be in the Fujiki class $\cC$.
\endclaim

By [Dem85], Conjecture (1.19) holds true in case $\alpha$ is an
integral class. Our hope is that the general case can be attained by
the diophantine approximation technique described earlier; there are 
however major hurdles, see [Lae02] for a few hints on these issues.

The author wishes to thank the organizers of the Riemann International
School of Mathematics held in Verbania in April 2009, for the opportunity of
publishing these notes in the RISM Proceedings volume.
\bigskip

\section{2. Variation of asymptotic cohomology groups}

We give here a proof a Theorem 1.8 in the context of general compact complex
manifolds~$X$. All norms occurring below are computed with respect to a fixed 
hermitian metric $\omega$ on $X$.

\claim (2.1) Lemma|Let $X$ be a compact complex $n$-fold.
Then for every coherent sheaf $\cF$ on X, there is a constant $C_{\cF}>0$
such that for every holomorphic line bundle $L$ on~$X$
we have
$$
h^q(X,\cF\otimes\cO_X(L))\le C_\cF(\Vert c_1(L)\Vert+1)^p
$$
where $p=\dim\Supp\cF$.
\endclaim

\proof. We prove the result by induction on $p\,$; it is indeed clear
for $p=0$ since we then have cohomology only in degree $0$ and the
dimension of $H^0(X,\cF\otimes\cO_X(L))$ does not depend on~$L$
when $\cF$ has finite support.
Let us consider the support $Y$ of $\cF$ and a resolution of singularity
$\mu:\widehat Y\to Y$ of the corresponding (reduced) analytic space. 
Then $\cF$ is an $\cO_Y$-module for some non
necessarily reduced complex structure $\cO_Y=\cO_X/\cJ$ on $J$. We can
look at the reduced structure $\cO_{Y,\red}=\cO_X/\cI$, $\cI=\sqrt{\cJ}$,
and filter $\cF$ by $\cI^k\cF$, $k\ge 0$. Since $\cI^k\cF/\cI^{k+1}\cF$ is a
coherent $\cO_{Y,\red}$-module, we can easily reduce the situation to the
case where $Y$ is reduced and $\cF$ is an $\cO_Y$-module. In that case
the cohomology $H^q(X,\cF\otimes\cO_X(L))=H^q(Y,\cF\otimes\cO_Y(L_{|Y}))$
just lives on the reduced space~$Y$.

Now, we have an injective sheaf morphism
$\cF\to\mu_\star\mu^*\cF$ whose cokernel $\cG$ has support in dimension~$<p$.
By induction on $p$, we conclude from the exact sequence that
$$
\big|h^q(X,\cF\otimes\cO_X(L))-h^q(X,\mu_\star\mu^*\cF\otimes\cO_X(L))\big|
\le C_1(\Vert c_1(L)\Vert+1)^{p-1}.
$$
The fonctorial morphisms
$$
\eqalign{
\mu^*&{}:H^q(Y,\cF\otimes\cO_Y(L_{|Y}))\to H^q(\widehat Y,\mu^\star\cF
\otimes\cO_{\widehat Y}(\mu^*L)_{|Y}),\cr
\mu_*&{}:H^q(\widehat Y,\mu^\star\cF
\otimes\cO_{\widehat Y}(\mu^*L)_{|Y})
\to H^q(Y,\mu_*\mu^\star\cF\otimes\cO_Y(L_{|Y}))\cr}
$$
yield a composition
$$\mu_*\circ\mu^*:H^q(Y,\cF\otimes\cO_Y(L_{|Y}))\to
H^q(Y,\mu_*\mu^\star\cF\otimes\cO_Y(L_{|Y}))
$$
induced by the natural injection $\cF\to\mu_\star\mu^*\cF$. This implies
$$
h^q(Y,\cF\otimes\cO_Y(L_{|Y}))\le 
h^q(\widehat Y,\mu^\star\cF\otimes\cO_{\widehat Y}(\mu^*L_{|Y}))
+C_1(\Vert c_1(L)\Vert+1)^{p-1}.
$$
By taking a suitable modification $\mu':Y'\to Y$ of the
desingularization $\widehat Y$, we can assume that $(\mu')^*\cF$ is
locally free modulo torsion. Then we are reduced to the case where
$\cF'=(\mu')^*\cF$ is a locally free sheaf on a smooth manifold $Y'$,
and $L'=(\mu')^*L_{|Y}$.
In this case, we apply standard analysis (e.g.\ [Dem85]) to conclude 
that $h^q(Y',\cF'\otimes\cO_{Y'}(L'))\le C_2(\Vert c_1(L')\Vert+1)^p$.
Since $\Vert c_1(L')\Vert\le C_3\Vert c_1(L)\Vert$ by pulling-back,
the statement follows easily.\qed

\claim (2.2) Corollary|For every irreducible divisor $D$ on $X$, there
exists a constant $C_D$ such that
$$
h^q(D,\cO_D(L_{|D}))\le C_D(\Vert c_1(L)\Vert+1)^{n-1}
$$
\endclaim

\proof. It is enough to apply Lemma~2.1 with $\cF=(i_D)_*\cO_D$ where
$i_D:D\to X$ is the injection.\qed

\claim (2.3) Remark|{\rm It is very likely that one can get an ``elementary''
proof of Lemma~2.1 without invoking resolutions of singularities, e.g.\ 
by combining the Cartan-Serre finiteness argument along with the standard
Serre-Siegel proof based ultimately on the Schwarz lemma. In this 
context, one would invoke $L^2$ estimates to get explicit bounds for 
the homotopy
operators between \v{C}ech complexes relative to two coverings 
$\cU=(B(x_j,r_j))$, $\cU'=(B(x_j,r_j/2))$ of $X$ by concentric balls.
By exercising enough care in the estimates, it is likely that one could
reach an explicit dependence $C_D\le C'\Vert D\Vert$ for the constant
$C_D$ of Corollary~2.2. The proof would of course become much more
technical than the rather naive brute force approach we have used.}
\endclaim

\medskip\noindent
{\bf (2.4) Proof of Theorem 1.8.}

\noindent
(a) We want to compare the cohomology of $L$ and $L'=L\otimes\cO(D)$ on $X$.
For this we write $D=D_+-D_-$, and compare the cohomology of the pairs
$L$ and $L_1=L\otimes\cO(-D_-)$ one one hand, and of $L'$ and 
$L_1=L'\otimes\cO(-D_+)$ on the other hand. Since
$\Vert c_1(\cO(D))\Vert\le C\Vert D\Vert$ by elementary potential theory,
we see that is is enough to consider the case of a negative divisor,
i.e.\ $L'=L\otimes\cO(-D)$, $D\ge 0$. If~$D$ is an irreducible
divisor, we use the exact sequence
$$
0\to L\otimes\cO(-D) \to L \to \cO_D\otimes L_{|D}\to 0
$$
and conclude by Corollary 2.2 that
$$
\eqalign{
\big|h^q(X,L\otimes\cO(-D))-h^q(X,L)\big|&\le 
h^q(D,\cO_D\otimes L_{|D})+h^{q-1}(D,\cO_D\otimes L_{|D})\cr
&\le 2C_D(\Vert c_1(L)\Vert+1)^{n-1}.\cr}
$$
For $D=\sum p_jD_j\ge 0$, we easily get by induction
$$
\big|h^q(X,L\otimes\cO(-D))-h^q(X,L)\big|\le 
2\sum_j p_jC_{D_j}\Big(\Vert c_1(L)\Vert+
\sum_k p_k\Vert\nabla_k\Vert+1\Big)^{n-1}
$$
If we knew that $C_D\le C'\Vert D\Vert$ as expected in Remark~2.3, then 
the argument would be complete without any restriction on $D$. 
The trouble disappears if we fix $D$ in
a finitely generated subgroup $\Gamma$ of divisors, because only finitely
many irreducible components appear in that case, and so we have to deal
with only finitely many constants $C_{D_j}$. Property (1.8~a) is proved.

\medskip\noindent
(b) Fix once for all a finite set of divisors $(\Delta_j)_{1\le j\le t}$
providing a basis of $\DNS_\bR(X)\subset H^{1,1}_\BC(X,\bR)$.
Take two elements $\alpha$ and $\beta$ in $\DNS_\bR(X)$, and fix
$\varepsilon>0$. Then $\beta-\alpha$ can be $\varepsilon$-approximated
by a $\bQ$-divisor $\sum\lambda_jD_j$, $\lambda_j\in \bQ$,
and we can find a pair $(k,L)$ with $k$ arbitrary large
such that ${1\over k}c_1(L)$ is $\varepsilon$-close to $\alpha$ and
${n!/k^n}h^q(X,L)$ approaches $\widehat h^q(X,\alpha)$ by $\varepsilon$.
Then ${1\over k}L+\sum \lambda_j\Delta_j$ approaches $\beta$ as closely
as we want. When approximating $\beta-\alpha$, we can arrange that 
$k\lambda_j$ is an integer by taking $k$ large enough. Then 
$\beta$ is approximated by ${1\over k}c_1(L')$ with
$L'=L\otimes\cO(\sum k\lambda_j\Delta_j)$. Property (a) implies
$$
\eqalign{
h^q(X,L')-h^q(X,L)
&\ge 
-C\Big(\Vert c_1(L)\Vert+\Big\Vert\sum k\lambda_j\Delta_j\Big\Vert\Big)^{n-1}
\Big\Vert\sum k\lambda_j\Delta_j\Big\Vert\cr
&\ge -Ck^n\big(\Vert\alpha\Vert+\varepsilon+\Vert\beta-\alpha\Vert+\varepsilon)^{n-1}
(\Vert\beta-\alpha\Vert+\varepsilon).\cr}
$$
We multiply the previous inequality by $n!/k^n$ and get in this way
$$
{n!\over k^n}h^q(X,L')\ge \widehat h^q(X,\alpha)-\varepsilon
 -C'\big(\Vert\alpha\Vert+\Vert\beta\Vert+\varepsilon)^{n-1}
(\Vert\beta-\alpha\Vert+\varepsilon).
$$
By taking the limsup and letting $\varepsilon\to 0$, we finally obtain
$$
\widehat h^q(X,\beta)-\widehat h^q(X,\alpha)\ge
-C'\big(\Vert\alpha\Vert+\Vert\beta\Vert)^{n-1}\Vert\beta-\alpha\Vert.
$$
Property (1.8~b) follows by exchanging the roles of $\alpha$ and $\beta$.\qed

\section{3. Monge-Amp\`ere volume formula}

The main goal of this section is to address the volume formula problem, namely
whether
$$
\Vol(\alpha)=\inf_{u\in\alpha}\int_{X(u,0)}u^n\qquad
\hbox{($u$ smooth)}\leqno(3.1)
$$
for every class $\alpha\in H^{1,1}_\BC(X,\bR)$ on a compact complex
$n$-fold~$(X,\omega)$.

\medskip\noindent
{\bf (3.2) Proof of the inequality${}\le$ (without restrictions)}

If $X$ does not admit any K\"ahler current, then the volume of every
class $\alpha$ is~$0$ and the inequality is trivially true. Therefore
we can assume that $X$ is in the Fujiki class~$\cC$. Then there exists a
K\"ahler modification $\mu:\widetilde X\to X$. Assume that we have a proof
for the K\"ahler case. Then
$$
\Vol(\alpha)=\Vol(\mu^*\alpha)\le\inf_{v\in\mu^*\alpha}\int_{\widetilde X(v,0)}v^n
\le \inf_{u\in\alpha}\int_{X(u,0)}u^n
$$
by restricting the inf to $v=\mu^*u$. This shows that it is enough to 
consider the case when $X$ is K\"ahler. We have something to prove
only when $\alpha\in\cE^+_X$, i.e.\ when $\alpha$ contains a K\"ahler 
current (a so-called ``big class''). Fix a $(1,1)$-form $u\in\alpha$.
We can then introduce
$$
\varphi(x):=\sup\big\{\psi(x)\,;\;\psi\le 0~\hbox{and}~u+i\ddbar\psi\ge 0
~\hbox{on $X$}\big\},\leqno(3.3)
$$
where the supremum is taken over all quasi-psh functions $\psi$
satisfying the given conditions $\psi\le 0$ and $u+i\ddbar\psi\ge 0$.
The following properties have been proved in [BD09] (cf.\ 
Theorem~1.4 and Corollary 2.5).

\claim (3.4) Lemma|Let $Z_0$ be the analytic set of poles of any K\"ahler
current $T_0\in\alpha$. Then $T=u+i\ddbar\varphi\ge 0$ and $\varphi$
is continuous with locally bounded second derivatives $\partial^2/\partial z_j
\partial\overline z_k$ on $X\ssm Z_0$. Moreover, if $S$ is the set of
points $z\in X\ssm Z_0$ where~$\varphi(z)=0$, then $S\subset\{z\,;u(z)\ge 0\}$
and
$$
\Vol(\alpha)=\int_Su^n=\int_{X\ssm Z_0}(u+i\ddbar\varphi)^n.
$$
\endclaim

\noindent
Since $S\subset \{z\,;u(z)\ge 0\}$, we 
immediately conclude from these equalities that
$$
\Vol(\alpha)\le \int_{\{z\,;u(z)\ge 0\}}u^n=
\int_{\{z\,;u(z)>0\}}u^n=\int_{X(u,0)}u^n.
$$

\noindent
{\bf (3.5) Proof of the volume formula for $\alpha\in\NS_\bR(X)$
on a projective surface.}

By definition, the volume $\Vol(\alpha)$ is obtained as the supremum
of $\int_{X\ssm \Sing(T)} T^n$ for K\"ahler currents with analytic
singularities in $\alpha$. By [Dem92] and [BDPP04], there exists a
blow-up $\mu:\widetilde X\to X$ such that $\mu^*T=[E]+\beta$ where $E$ is a
normal crossing divisor on $\widetilde X$ and $\beta\ge 0$
smooth. Until now, this is valid for an arbitrary compact complex
manifold~$X$. If moreover $X$ is projective and $\alpha\in\NS_\bR(X)$,
it is shown in [BDPP04] that we have the ``orthogonality property''
$$
[E]\cdot \beta^{n-1}=\int_{E}\beta^{n-1}\le 
C\big(\Vol(\alpha)-\beta^n\big)^{1/2},\leqno(3.6)
$$
in other words, $E$ and $\beta$ become ``more and more orthogonal'' as
$\beta^n$ approaches the volume. Our method consists of approaching
$[E]+\beta$ by smooth closed $(1,1)$-forms $u_\varepsilon$ in the same 
$\ddbar$-cohomology class as $[E]+\beta$, in such a way that 
$$\int_{\widetilde X(u_\varepsilon,0)}u_\varepsilon^n$$
will not be substantially larger than the volume $\int_{\widetilde X}\beta^n$.
For this, we select a hermitian metric $h$ on $\cO(E)$ and put
$$
u_\varepsilon={i\over 2\pi}\ddbar\log(|\sigma_E|_h^2+\varepsilon^2)+
\Theta_{\cO(E),h}+\beta\leqno(3.7)
$$
where $\sigma_E\in H^0(\widetilde X,\cO(E))$ is the canonical section
and $\Theta_{\cO(E),h}$ the Chern curvature form. Clearly,by the
Lelong-Poincar\'e equation, $u_\varepsilon$ converges to $[E]+\beta$ in 
the weak topology as $\varepsilon\to 0$. Straightforward calculations yield
$$
u_\varepsilon={i\over 2\pi}{\varepsilon^2D^{1,0}_h\sigma_E\wedge 
\overline{D^{1,0}_h\sigma_E}\over(\varepsilon^2+|\sigma_E|^2)^2}
+{\varepsilon^2\over\varepsilon^2+|\sigma_E|^2}\Theta_{E,h}+\beta.
$$
The first term converges to $[E]$ in the weak topology, while the second,
which is close to $\Theta_{E,h}$ near~$E$, converges pointwise everywhere to $0$
on $\smash{\widetilde X}\ssm E$. A simple asymptotic analysis shows that 
$$
\Big({i\over 2\pi}{\varepsilon^2D^{1,0}_h\sigma_E\wedge 
\overline{D^{1,0}_h\sigma_E}\over(\varepsilon^2+|\sigma_E|^2)^2}
+{\varepsilon^2\over\varepsilon^2+|\sigma_E|^2}\Theta_{E,h}\Big)^p\to
[E]\wedge\Theta_{E,h}^{p-1}
$$
in the weak topology for $p\ge 1$, hence
$$
\lim_{\varepsilon\to 0}u_\varepsilon^n=
\beta^n+\sum_{p=1}^n{n\choose p}[E]\wedge \Theta_{E,h}^{p-1}
\wedge\beta^{n-p}.\leqno(3.8)
$$
In arbitrary dimension, the signature of $u_\varepsilon$ is hard to evaluate,
and it is also non trivial to decide the sign of the limiting measure
$\lim u_\varepsilon^n$. However, when $n=2$, we get the simpler formula
$$
\lim_{\varepsilon\to 0}u_\varepsilon^2=\beta^2+2[E]\wedge\beta
+[E]\wedge\Theta_{E,h}.
$$
In this case, $E$ can be assumed to be an exceptional divisor (otherwise
some part of it would be nef and could be removed from the poles of $T$).
Hence the matrix $(E_j\cdot E_k)$ is negative definite and we can find 
a hermitian
metric $h$ on $\cO(E)$ such that $(\Theta_{E,h})_{|E}<0$. Then 
$[E]\wedge\Theta_{E,h}$, which is the limit of the product of the first
two terms in $u_\varepsilon^2$, contributes negatively to the limit; all
other terms are nonnegative or have a mass converging to~$0$.
From this, one can easily infer by (3.6) that
$$
\limsup_{\varepsilon\to 0}\int_{\widetilde X(u_\varepsilon,0)}u_\varepsilon^2
\le \int_{\widetilde X}\beta^2+2[E]\wedge\beta
\le \Vol(\alpha)+2C(\Vol(\alpha)-\beta^2)^{1/2}.
$$
This is arbitrary close to $\Vol(\alpha)$ when $\beta^2$ approaches
the volume, and so property (1.18~b) is proved in dimension~$2$. Obviously
the $n$-dimensional case would require a deeper analysis of
``higher order'' orthogonality relations.\qed

\section{References}
\medskip

{\eightpoint

\bibitem[AG62]&Andreotti, A., Grauert, H.:& Th\'eor\`emes de finitude
  pour la cohomologie des espaces complexes;& Bull.\ Soc.\ Math.\
  France {\bf 90} (1962) 193--259&

\bibitem[BL04]&Birkenhake, Ch., Lange, H.:& Complex Abelian Varieties;&
  Second augmented edition, Grundlehren der Math.\ Wissenschaften,
  Springer, Heidelberg, 2004&

\bibitem[Bou02]&Boucksom, S.:& On the volume of a line bundle;&
  Internat.\ J.\ Math.\ {\bf 13} (2002), 1043--1063&

\bibitem[BD09]&Berman, R., Demailly, J.-P.:& Regularity of
  plurisubharmonic upper envelopes in big cohomology classes;&
  arXiv: math.CV/0905.1246, to appear in the Proceedings of
  the volume ``Geometry and Topology'' in honor of Oleg Y.\ Viro, edited 
  by B.\ Juhl-J\"oricke and M.\ Passare, Birkha\"user&

\bibitem[BDPP04]&Boucksom, S., Demailly, J.-P., P\u{a}un, M.,
  Peternell, Th.:& The pseudo-effective cone of a compact K\"ahler
  manifold and varieties of negative Kodaira dimension;& arXiv: 
  math.AG/0405285, see also Proceedings of the ICM 2006 in Madrid&

\bibitem[Dem85]&Demailly, J.-P.:& Champs magn\'etiques et
  in\'egalit\'es de Morse pour la $d''$-cohomologie;& Ann.\ Inst.\
  Fourier (Grenoble) {\bf 35} (1985), 189--229&

\bibitem[Dem91]&Demailly, J.-P.:& Holomorphic Morse inequalities; &
  Lectures given at the AMS Summer Institute on Complex Analysis held
  in Santa Cruz, July 1989, Proceedings of Symposia in Pure
  Mathematics, Vol.~{\bf 52}, Part~2 (1991), 93--114&

\bibitem[Dem90]&Demailly, J.-P.:& Singular hermitian metrics on
  positive line bundles;& Proceedings of the Bayreuth conference
  ``Complex algebraic varieties'', April~2-6, 1990, edited by
  K.~Hulek, T.~Peternell, M.~Schneider, F.~Schreyer, Lecture Notes in
  Math.\ n${}^\circ\,$1507, Springer-Verlag, 1992&

\bibitem[Dem92]&Demailly, J.-P.:& Regularization of closed positive currents 
  and Intersection Theory;& J.\ Alg.\ Geom.\ {\bf 1} (1992), 361--409&

\bibitem[DEL00]&Demailly, J.-P., Ein, L., {\rm and} Lazarsfeld, R.:& A
  subadditivity property of multiplier ideals;& Michigan Math. J. 48
  (2000), 137\--156&

\bibitem[DP04]&Demailly, J.-P., P\u{a}un, M:& Numerical
  characterization of the K\"ahler cone of a compact K\"ahler
  manifold;& arXiv: math.AG/0105176$\,$; Annals of Math.\ {\bf 159} (2004)
  1247--1274&

\bibitem[FKL07]&de Fernex, T., K\"uronya, A., Lazarsfeld, R.:& Higher
  cohomology of divisors on a projective variety;& Math.\ Ann.\ {\bf
  337} (2007) 443--455&

\bibitem[Fuj94]&Fujita, T.: & Approximating Zariski decomposition of
  big line bundles;& Kodai Math.\ J.\ {\bf 17} (1994) 1-\-3&

\bibitem[Hir54]&Hirzebruch, F.:& Arithmetic genera and the theorem of
  Riemann-Roch for algebraic varieties;& Proc.\ Nat.\ Acad.\ Sci.\
  U.S.A.\ {\bf 40} (1954), 110--114&

\bibitem[Hir56]&Hirzebruch, F.:& Neue topologische Methoden in der
  algebraischen Geometrie;& Ergebnisse der Mathematik und ihrer
  Grenzgebiete (N.F.), Heft 9. Springer-Verlag,
  Berlin-G\"ottingen-Heidelberg, 1956.165 pp; English translation: {\it
 Topological methods in algebraic geometry}; Springer-Verlag,
  Berlin (1966)&

\bibitem[Kur06]&K\"uronya, A.:& Asymptotic cohomological functions on
  projective varieties;& Amer.\ J.\ Math.\ {\bf 128} (2006) 1475--1519&

\bibitem[Lae02]&Laeng, L.:& Estimations spectrales asymptotiques en
  g\'eom\'etrie hermitienne;& Th\`ese de Doctorat de l'Universit\'e 
  de Grenoble~I, octobre 2002,\\
  http://www-fourier.ujf-grenoble.fr/THESE/ps/laeng.ps.gz and\\
  http://tel.archives-ouvertes.fr/tel-00002098/en/&

\bibitem[Laz04]&Lazarsfeld, R.:& Positivity in Algebraic Geometry
  I.-II;& Ergebnisse der Mathematik und ihrer Grenzgebiete,
  Vols. 48-49., Springer Verlag, Berlin, 2004&

\bibitem[Rie57]&Riemann, B.:& Theorie der Abel'schen Functionen;& J.\
  f\"ur Math.\ {\bf 54} (1857), in Gesam\-melte mathematische Werke
  (1990) 120--144&

\bibitem[Roc65]&Roch, G.:& \"Uber die Anzahl der willkurlichen
  Constanten in algebraischen Functionen;& J.\ f\"ur Math.\ {\bf 64} 
  (1865) 372--376&

}
\vskip20pt
\noindent
(version of February 23, 2010, printed on \today)
\end